\documentclass{amsproc}

\usepackage{a4wide}
\usepackage{amsmath}
\usepackage{enumerate}

\usepackage{amsmath,amsthm}
\usepackage{amsfonts}
\usepackage{amssymb}
\usepackage{longtable}
\usepackage[matrix,arrow,curve]{xy}

\usepackage{amsmath}
\usepackage{amsfonts}
\usepackage{amssymb}

\newtheorem{theorem}{Theorem}
\newtheorem{cor}{Corollary}

\newtheorem{ex}{Example}{\bf}{\rm}

{\bf}{\rm}{\rm}

\newtheorem*{theorem*}{Theorem}
\newtheorem*{cor*}{Corollary}

\def\Real{\mathbb{R}}
\def\SO{\text{\rm SO}}
\def\Spin{\text{\rm Spin}}
\def\CO{\text{\rm CO}}
\def\CSpin{\text{\rm CSpin}}
\def\Hol{\text{\rm Hol}}

\newcommand{\frakg}{\mathfrak{g}}
\newcommand{\frakh}{\mathfrak{h}}
\newcommand{\mathR}{\mathbb{R}}
\newcommand{\frakso}{\mathfrak{so}}
\newcommand{\frakco}{\mathfrak{co}}

\newcommand{\g}{\mathfrak{g}}
\newcommand{\h}{\mathfrak{h}}
\newcommand{\so}{\mathfrak{so}}
\newcommand{\co}{\mathfrak{co}}
\newcommand{\hol}{\mathfrak{hol}}

\newcommand{\zr}{\ltimes}

\def\id{\mathop\text{\rm id}\nolimits}

\def\Ric{\mathop\text{\rm Ric}\nolimits}

\usepackage{stackrel}

\newcommand{\be}{\begin{equation}}
\newcommand{\ee}{\end{equation}}

\let\ge=\geqslant
\let\leq=\leqslant
\let\geq=\geqslant

\begin{document}

\title{Parallel spinors on Lorentzian Weyl spaces}

\author{Andrei Dikarev}\thanks{$^1$ Masaryk University, Faculty of Science, Department of Mathematics and Statistics, Kotl\'a\v{r}sk\'a 2, 611 37 Brno, Czech
	Republic}

\author{Anton S. Galaev}\thanks{$^2$University of Hradec Kr\'alov\'e, Faculty of Science, Rokitansk\'eho 62, 500~03 Hradec Kr\'alov\'e,  Czech
Republic\\
E-mail: anton.galaev(at)uhk.cz}

\begin{abstract}

Holonomy algebras of Lorentzian Weyl spin manifolds with weighted parallel spinors are found. 
For Lorentzian Weyl manifolds admitting recurrent null vector fields are introduced special local coordinates similar to Kundt and Walker ones. Using that, the local form of all  Lorentzian Weyl spin manifolds with weighted parallel spinors is given. The Einstein-Weyl equation for the obtained Weyl structures is analyzed and examples of Einstein-Weyl spaces with weighted parallel spinors are given.

\vskip0.5cm


\vskip0.5cm



\end{abstract}

\maketitle

\section*{Introduction} The present paper is dedicated to the study of weighted parallel spinors on Weyl spin manifolds of Lorentzian signature. 
Parallel spinors are special Killing spinors which represent supersymmetry generators  of supersymmetric field theories and supergravity theories. The physical motivation for study Weyl spaces with weighted parallel spinors may be found in \cite{MOP}.

There is a one-to-one correspondence between the parallel spinors and the holonomy-invariant elements of the spinor module. This correspondence allowed to describe the following  simply connected spin manifolds with parallel spinors: Riemannian manifolds \cite{Hitchin,Wang}, pseudo-Riemannian manifolds with irreducible holonomy groups \cite{Baum-Kath}, Lorentzian manifolds \cite{Br,FF,L01,L07}.   

The work \cite{Moroianu} initiated the study of parallel spinors on Weyl manifolds. It was  proven that in the Riemannian signature,  parallel spinors of weight zero for non-closed Weyl structures exist only on certain non-compact 4-dimensional manifolds. For the exposition of this and many other related results see the recent book \cite{Spinor}. Killing spinors of arbitrary weight on Riemannian  Weyl spin manifolds were studied for the first time in \cite{Buchholz}. 

The work  \cite{MOP} provides a deep investigation of weighted parallel spinors on Einstein-Weyl manifolds of Lorentzian signature with a special attention to the dimensions 4 and 6. The techniques developed for classification of supergravity solutions was used in that work.

In the present paper we provide a description of simply connected Lorentzian Weyl spin manifolds admitting  weighted parallel spinors. The main tool for that are holonomy groups. In Sections \ref{sec1} and \ref{sec1.1} we give necessary background on weighted parallel spinors and spinor modules. In Section \ref{sec2} we review the recent classification of holonomy algebras of Lorentzian Weyl manifolds~\cite{Andr}. In Section \ref{sec3} we classify the holonomy algebras of Lorentzian Weyl spaces admitting weighted parallel spinors. It turns out that for non-closed Weyl structures, there are  two types of such algebras. In each case, the dimension of the space of parallel spinors is found. In Section \ref{sec4} we give a local description of Lorentzian Weyl structures that admit parallel distributions of null lines.   We show that under a mild condition, such structures may be described by Walker coordinates~\cite{Walkerbook}. Such local coordinates were considered in particular cases in, e.g., \cite{DMT,MOP}.
In Section \ref{sec5}, we found the local form of the Lorentzian Weyl structures admitting weighted parallel  spinors. 
In Section \ref{secEW}, we analyze the Einstein-Weyl equation for the obtained Weyl structures  and construct examples of Einstein-Weyl spaces with weighted parallel spinors. Some examples have previously appeared in \cite{DMT,MOP} and other literature. It turns out that the Einstein-Weyl equation implies that the weight of a non-zero weighted parallel spinor is equal to $\dim M-4$. 
Parallel spinors of that weight were studied in \cite{MOP}. In contrast, Section~\ref{sec5} provides Weyl structures with non-zero weighted  parallel spinors of arbitrary weight. 

As it is explained in \cite{Ag-Fr}, the results on Weyl connections may be applied to metric connections with vectorial torsion.

\section{Weighted parallel spinors}\label{sec1}

Let $(M,c,\nabla)$ be a Weyl manifold, i.e., $c$ is a conformal class of pseudo-Riemannian metrics of signature $(r,s)$ ($r$ denotes the number of minuses) on a smooth manifold $M$, and $\nabla$ is a torsion-free affine connection on $M$ such that, for each $g\in c$, there exists a 1-form $\omega$ with \begin{equation}\label{nablag}\nabla g=2\omega\otimes g.\end{equation}
Note that a metric $g$ and the corresponding 1-form $\omega$ determine the connection $\nabla$, it holds \begin{equation}\label{formulaK}\nabla=\nabla^g+K,\quad 	g(K_X(Y),Z)=g(Y,Z)\omega(X)+g(X,Z)\omega(Y)-g(X,Y)\omega(Z),\end{equation} where $\nabla^g$ is the Levi-Civita connection of the metric $g$, and $X,Y,Z$ are vector fields on $M$.

We  follow the exposition from \cite{Buchholz} addapting it to the pseudo-Riemannian settings. Let $$\CSpin(r,s)=\Real_+\times \Spin(r,s),\quad   \CO(r,s)=\Real_+\times \SO(r,s).$$ We obtain the twofold covering
$$\lambda^c:\CSpin(r,s)\to \CO(r,s).$$ 
For an oriented Weyl manifold $(M,c,\nabla)$, let $P_\CO$ be the corresponding bundle of oriented conformal  frames.
A spin structure on $(M,c,\nabla)$ is a $\lambda^c$-reduction of the bundle $P_\CO$ to a bundle $P_\CSpin$. The existence of a spin structure on  
$(M,c,\nabla)$ is equivalent to the existence of a spin structure on $(M,g)$ for any $g\in c$, since the maximal compact subgroups of   $\CSpin(r,s)$ and $\Spin(r,s)$ coincide. The conditions for a pseudo-Riemannian manifold to admit a spin structure may be found in \cite{Baum}, see also \cite{Al-Ch}. 

Let $\Delta_{r,s}$ be the standard spinor module, we recall its definition in the next section. Let $w\in\Real$. The spinor representation $k^w$ of $\CSpin(r,s)$ on $\Delta_{r,s}$ with weight $w$ is defined by $$k^w(a,h)\psi=a^wh\psi,$$
where $a\in\Real_+$, $h\in \Spin(r,s)$, and $\psi\in\Delta_{r,s}$.

The spinor bundle with weight $w$ is defined by $$S^w=P_\CSpin\times_{k^w}\Delta_{r,s}.$$
The connection $\nabla$ induces a connection $\nabla^S$ on $S^w$. 
Let us fix a metric $g\in c$. Denote by $\nabla^g$ and $\nabla^{S,g}$ the corresponding Levi-Civita connection and the connection in the spinor bundle, respectively. The connection $\nabla^S$ is given by
$$\nabla_X^S\psi=\nabla^{S,g}_X\psi-\frac{1}{2}X\cdot\omega\cdot\psi+\left(w-\frac{1}{2}\right) \omega(X)\psi.$$
A spinor  $\psi$ of weight $w$ on $M$ is called parallel if $\nabla^S\psi=0$.

For the holonomy group of the connection $\nabla^S$ at a point $x\in M$  it holds $$\Hol_x(\nabla^S)=(\lambda^c)^{-1}(\Hol_x(\nabla)),$$ where $\Hol_x(\nabla)$ is the holonomy group of the connection $\nabla$ at $x\in M$.
 Clearly, there is an isomorphism of  the  space of parallel spinors of weight $w$ and the space of holonomy-invariant spinors of the spinor module of weight $w$:
$$\{\psi\in\Gamma(S^w)|\nabla^S\psi=0\}\cong\{\psi_x\in S_x=\Delta_{r,s}|A\psi_x=\psi_x\,\,\forall A\in \Hol_x(\nabla^S)\}.$$ If the manifold $M$ is simply connected, then the both spaces are isomorphic to 
$$\{\psi_x\in S_x=\Delta_{r,s}|\xi\psi_x=0\,\,\forall \xi\in \hol_x(\nabla^S)\},$$ where $\hol_x(\nabla^S)$ is the holonomy algebra of the connection $\nabla^S$ at the point $x$, i.e., the Lie algebra of the holonomy group $\Hol_x(\nabla^S)$.

\section{Spinor modules}\label{sec1.1}

Let us fix some standard notation.
Let $\mathbb{R}^{r,s}$ be a pseudo-Euclidean space with the metric
$g$ of signature $(r,s)$.  Let
$({\mathcal C}l_{r,s},\cdot)$ be the corresponding Clifford algebra
and $\mathbb{C} l_{r,s}= {\mathcal C}l_{r,s}\otimes\mathbb{C}$ be its
complexification. The last algebra can be represented as a matrix
algebra in the following way. Consider the basis
$$\left(u(\epsilon)=\frac{\sqrt{2}}{2}\left(\begin{array}{c}1\\-\epsilon{\rm i}\end{array}\right),\,
\epsilon=\pm 1\right)$$ of $\mathbb{C}^2$. Define the following
isomorphisms of $\mathbb{C}^2$:
$$E={\rm id},\quad T=\left(\begin{array}{cc}0&-{\rm i}\\{\rm i}&0\end{array}\right),
\quad U=\left(\begin{array}{cc}{\rm i}&0\\0&-{\rm
	i}\end{array}\right), \quad V=\left(\begin{array}{cc}0&{\rm
	i}\\{\rm i}&0\end{array}\right).$$ It holds
$$T^2=-V^2=-U^2=E,\quad UT=-{\rm i} V,\quad VT={\rm i} U,\quad UV=-{\rm i} T,$$
$$Tu(\epsilon)=-\epsilon u(\epsilon),\quad Uu(\epsilon)={\rm i}
u(-\epsilon),\quad Tu(\epsilon)=\epsilon u(-\epsilon).$$ Let
$n=r+s$. A basis $e_1,...,e_n$ of $\mathbb{R}^{r,s}$ is called
orthonormal if $g(e_i,e_j)=k_i\delta_{ij}$, where $k_i=-1$ if
$1\leq i\leq r$, and $k_i=1$ if $r+1\leq i\leq n$. Let us fix such
basis. For an integer $m$ denote by $\mathbb{C}(m)$ the algebra of
the complex square matrices of order $m$. Define the following
isomorphisms:
\begin{itemize}
	\item[1)] if $n$ is even, then define $\Phi_{r,s}:\mathbb{C} l_{r,s}\to\mathbb{C}\Big(2^{\frac{n}{2}}\Big)$ by
	\begin{align}\Phi_{r,s}(e_{2k-1})&=\tau_{2k-1}E\otimes\cdots\otimes E\otimes U\otimes
	\underbrace{T\otimes \cdots\otimes T}_{(k-1)-times},\label{Phi1}\\
	\Phi_{r,s}(e_{2k})&=\tau_{2k}E\otimes\cdots\otimes E\otimes
	V\otimes \underbrace{T\otimes \cdots\otimes
		T}_{(k-1)-times},\label{Phi2}\end{align} where $1\leq k\leq
	\frac{n}{2}$,  $\tau_i={\rm i}$ if $1\leq i\leq r$, and $\tau_i=1$
	if $r+1\leq i\leq n$;
	\item[2)] if $n$ is odd, then define $\Phi_{r,s}:\mathbb{C} l_{r,s}\to\mathbb{C}\Big(2^{\frac{n-1}{2}}\Big)\oplus\mathbb{C}\Big(2^{\frac{n-1}{2}}\Big)$ by
	\begin{align}\Phi_{r,s}(e_{k})&=(\Phi_{r,s-1}(e_{k}),\Phi_{r,s-1}(e_{k})),\quad k=1,...,n-1,\\
	\Phi_{r,s}(e_{n})&=({\rm i}T\otimes\cdots\otimes T,-{\rm i}
	T\otimes\cdots\otimes T).
	\end{align}
\end{itemize}
The obtained representation space
$\Delta_{r,s}=\mathbb{C}^{2^{[\frac{n}{2}]}}$ is called {\it the
	spinor module}.
We write $A\cdot s=\Phi_{r,s}(A)s$ for all
$A\in\mathbb{C} l_{r,s}$, $s\in\Delta_{r,s}$. We will consider the
following basis of  $\Delta_{r,s}$:
$$(u(\epsilon_k,...,\epsilon_1)=u(\epsilon_k)\otimes\cdots\otimes u(\epsilon_1)|\epsilon_i=\pm 1).$$
Recall that the Lie algebra $\mathfrak{spin}(r,s)$ of the Lie
group ${\rm Spin}(r,s)$ can
be embedded into $\mathbb{C} l_{r,s}$ in the following way:
$$\mathfrak{spin}(r,s)={\rm span}\{e_i\cdot e_j|1\leq i<j\leq n\}.$$
The Lie algebra $\mathfrak{so}(r,s)$ can be identified with the
space of bivectors $\Lambda^2\mathbb{R}^{r,s}$ in such a way that
\begin{equation}\label{isom_so}(x\wedge y)z=g(x,z)y-g(y,z)x,\quad x,y,z\in\mathbb{R}^{r,s}.\end{equation}
There is the isomorphism
$$\lambda_*:\mathfrak{so}(r,s)\to\mathfrak{spin}(r,s),\quad
\lambda_*(x\wedge y)=x\cdot y.$$ The obtained representation of
$\mathfrak{so}(r,s)$ in $\Delta_{r,s}$ is irreducible if $n$ is
odd, and this representation   splits into the direct sum of two
irreducible modules $$\Delta_{r,s}^\pm ={\rm
	span}\{u(\epsilon_k,...,\epsilon_1)|\epsilon_1=\cdots=\epsilon_k=\pm
1\}$$ if $n$ is even.

\section{Holonomy of Lorentzian Weyl manifolds}\label{sec2}

Let us summarize the recent result on the classification of the holonomy algebras of Weyl manifolds of Lorentzian signature \cite{Andr}. Let $(M,c,\nabla)$ be a Weyl manifold of Lorentzian signature $(1,n+1)$, $n\geq 1$. Then its holonomy algebra is contained in the conformal Lorentzian algebra $$\co(1,n+1)=\Real\id\oplus\so(1,n+1).$$ If the Weyl structure is closed, then the holonomy algebra is contained in $\so(1,n+1)$ and it is well-studied \cite{L07}. By that reason we suppose that the Weyl structure is non-closed, and the holonomy algebra is not contained in  $\so(1,n+1)$.

Fix a Witt basis $p,e_1,\dots,e_n,q$ of the Minkowski space $\Real^{1,n+1}$. With respect to that basis the subalgebra of $\so(1,n+1)$ preserving the null line $\Real p$ has the following matrix form:

\begin{equation*}
\frakso (1, n + 1)_{\mathR p} = \left\{
 \left. 
\begin{pmatrix}     
a & X^t & 0 \\
0 & A & -X \\
0 & 0 & -a
\end{pmatrix} \right|
\begin{matrix}     
a \in \mathR \\
A \in \frakso (n) \\
X \in \mathR^n
\end{matrix} \right\} .
\end{equation*}	
We get the decomposition 
$$\frakso (1, n + 1)_{\mathR p} = (\mathbb{R}\oplus\frakso(n))\ltimes\mathbb{R}^n.$$
An element of $$\frakco (1, n + 1)_{\mathR p}=\Real\id\oplus\frakso (1, n + 1)_{\mathR p}$$ will be denoted by $(b, a, A, X)$, where $b \in \mathR$ and  $(a, A, X)$ is defined by the above matrix.

The only irreducible holonomy algebra is the whole $\co(1,n+1)$.
The first result describes the holonomy of conformal products (in the sense of \cite{Moroianu1}).

\begin{theorem} \cite{Andr}
	\label{class2Th}
	If the holonomy algebra $\frakg\subset  \frakco (1, n + 1)$ of a non-closed Weyl connection preserves  an orthogonal decomposition  $$\mathR^{1, n + 1} = \mathR^{1, k + 1} \oplus \mathR^{n - k}, \quad -1 \leqslant k \leqslant n - 1,$$ then  $\frakg$ is one of the following:
	\begin{itemize}
		\item $\mathR \id\oplus \frakso (1, k + 1) \oplus \frakso (n - k) $, $ \quad -1 \leqslant k \leqslant n-1$;
		\item $ \mathR (1, -1, 0, 0) \oplus \mathfrak{k} \oplus \frakso (n - k) \ltimes \mathR^{k}\subset\frakco (1, n + 1)_{\mathR p}$,  $0\leq k\leq n-1$;
		\item $ \{ (b, a, 0, 0)|b,a\in\Real\} \oplus \mathfrak{k} \oplus \frakso (n - k) \ltimes \mathR^{k}\subset\frakco (1, n + 1)_{\mathR p}$,  $1\leq k\leq n-1$.
			\end{itemize}
Here $\mathfrak{k} \subset \frakso (k)$ is the holonomy algebra of a Riemannian manifold.
\end{theorem}

The rest of the holonomy algebras do not preserve any non-degenerate subspace, but preserve a null line and may be considered as  subalgebras of $\co(1,n+1)_{\Real p}$. 

\begin{theorem}
	\label{class3Th} \cite{Andr}
	If the holonomy algebra $\frakg\subset  \frakco (1, n + 1)$ of a non-closed Weyl connection does not preserve any non-degenerate subspace of $\Real^{1,n+1}$ and preserves a null line, then $\frakg$ is one of the following:
	\begin{itemize}
	\item $\mathfrak{g}^{\Real,1,\mathfrak{h}} = \{ (b, a, 0, 0) \mid b,a \in \mathR \}\oplus \h \ltimes \mathR^{n}$;
	
	\item $\mathfrak{g}^{\Real,2,\mathfrak{h}} = \{ (b, 0, 0, 0) \mid b \in \mathR \}\oplus \h \ltimes \mathR^{n}$;
	
	\item $\mathfrak{g}^{\Real,3,\mathfrak{h},\varphi} = \{ (b, \varphi (A), A, 0) \mid b\in\Real,\,A \in \frakh \} \ltimes \mathR^{n}$;

	\item $\mathfrak{g}^{\beta,\theta,1,\mathfrak{h}} = \{(\beta a+ \theta ( A), a, A, 0) \mid a \in \mathR, A \in \frakh \} \ltimes \mathR^{n}$;
	
		\item $\mathfrak{g}^{\theta,2,\mathfrak{h}} = \{ (\theta (A), 0, A, 0) \mid A \in \frakh \} \ltimes \mathR^{n}$;
		
		\item $\mathfrak{g}^{\theta,3,\mathfrak{h},\varphi} = \{ (\theta (A), \varphi (A), A, 0) \mid A \in \frakh \} \ltimes \mathR^{n}$.

			\end{itemize}
	Here  $\frakh\subset\so(n)$ is the holonomy algebra of a Riemannian manifold, $\beta\in\Real$, and
$\theta, \varphi : \frakh \rightarrow \mathR$  are  linear maps such that $\theta \big|_{[\frakh, \frakh]} =\varphi\big|_{[\mathfrak{h}, \mathfrak{h}]}=0$, $\varphi\neq 0$. Moreover, for the forth algebra it holds $\beta\neq 0$ or $\theta\neq 0$, and for the last two algebras it holds $\theta\neq 0$.
\end{theorem}

\section{Weighted parallel spinors and holonomy}\label{sec3}

In this section we characterize  simply connected Weyl manifolds   of Lorentzian signature admitting weighted parallel spinors in terms of the holonomy algebras.

First we recall some known results. Let $\h\subset\so(n)$ be the holonomy algebra of a Riemannian manifold. Recall that
there is an orthogonal decomposition
\begin{equation}\label{LM0A}\mathbb{R}^{n}=\mathbb{R}^{n_{0}}\oplus\mathbb{R}^{n_1}\oplus\cdots\oplus\mathbb{R}^{n_s}\end{equation} and the
corresponding decomposition into the direct sum of ideals
\begin{equation}\label{LM0B}\mathfrak{h}=\{0\}\oplus\mathfrak{h}_1\oplus\cdots\oplus\mathfrak{h}_s\end{equation} such that $\mathfrak{h}$ annihilates
$\mathbb{R}^{n_{0}}$, $\mathfrak{h}_i(\mathbb{R}^{n_j})=0$ for
$i\neq j$, and $\mathfrak{h}_i\subset\mathfrak{so}(n_i)$ is an
irreducible subalgebra for $1\leq i\leq s$. Moreover, the
subalgebras $\mathfrak{h}_i\subset\mathfrak{so}(n_i)$ are the
holonomy algebras of Riemannian manifolds. This decomposition of the holonomy algebra corresponds to the de~Rham decomposition of the manifold $(M,g)$.
Let $(M,g)$ be a simply connected  Riemannian spin manifold of dimension $n$ with the holonomy algebra $\mathfrak{h}\subset\mathfrak{so}(n)$. Then $(M,g)$ admits a non-zero parallel spinor if and only if 
in the decomposition \eqref{LM0B} of
$\mathfrak{h}\subset\mathfrak{so}(n)$ each subalgebra
$\mathfrak{h}_i\subset\mathfrak{so}(n_i)$ coincides with one of
the Lie algebras
$\mathfrak{su}(\frac{n_i}{2})$, $\mathfrak{sp}(\frac{n_i}{4})$, $G_2\subset\mathfrak{so}(7)$, $\mathfrak{spin}(7)\subset\mathfrak{so}(8)$ \cite{Wang}.

For the holonomy algebra $\g\subset\so(1,n+1)$ of a Lorentzian manifold $(M,g)$, the Wu Theorem implies  decompositions similar to \eqref{LM0A} and \eqref{LM0B} with the following difference. There are two possible cases. In the first case, the subspace $\Real^{n_0}$ annihilated by $\g$ is of Lorentzian signature, then $\g$ is just the same as the above $\h$, and the condition implied by the existence of a non-zero parallel tensor is the same as above. In the second case, one of the factors, e.g., $\Real^{n_1}$ is of Lorentzian signature, i.e., $\h_1\subset\so(1,n_1-1)$, and $\h_1$ does not preserve any proper non-degenerate vector subspace of $\Real^{1,n_1-1}$. There exists a non-zero parallel spinor   if and only if $\h_2,\dots,\h_s$ are as above and $\h_1$ annihilates a null vector in $\Real^{1,n_1-1}$, and it is isomorphic to the Lie algebra $\mathfrak{f}\zr\Real^{n_1-2}$, where $\mathfrak{f}\subset \so(n_1-2)$ is the holonomy algebra of a  Riemannian spin manifold admitting a non-zero parallel spinor ($n_1\geq 3$). These results are due to \cite{L01,L07}.

\vskip0.3cm

Now we prove one of the main results of the present paper.

\begin{theorem}\label{thHolSpinor}
Let  $(M,c,\nabla)$ be a  simply connected  manifold with a non-closed      Weyl spin  structure of Lorentzian signature $(1,n+1)$, $n\geq 1$.
Then it admits a non-zero parallel spinor of weight~$w$ if and only if  its holonomy algebra is one of the following subalgebras of $\co(1,n+1)_{\Real p}$:
\begin{itemize}
\item $\g^{w,\h}=\mathR \left(1, \frac{w}{2}, 0, 0\right) \oplus \mathfrak{h} \ltimes \mathR^{n} $ for arbitrary $w$;
	\item 
$\g^\mathfrak{k}=\mathR (1, -1, 0, 0) \oplus \mathfrak{k} \ltimes \mathR^{n-1} $ only for $w=-2$.
\end{itemize}
Here  $\mathfrak{k} \subset \frakso (n-1)$ and $\mathfrak{h} \subset \frakso (n)$ are the holonomy algebras of Riemannian spin manifolds  carrying parallel spinors. 

The dimension of the space of weighted parallel spinors on  $(M,c,\nabla)$ with the holonomy algebra $\g^{w,\h}$ equals to
$$\dim\{\psi\in\Delta_{n}|\xi\psi=0\,\,\forall\xi\in\h\}.$$
The dimension of the space of  weighted parallel spinors on  $(M,c,\nabla)$ with the holonomy algebra $\g^{\mathfrak{k}}$ equals to
$$\left\{\begin{matrix}\dim\{\psi\in\Delta_{n-1}|\xi\psi=0\,\,\forall\xi\in\mathfrak{k}\},\quad \text{if $n$  is even},\\
	2\dim\{\psi\in\Delta_{n-1}|\xi\psi=0\,\,\forall\xi\in\mathfrak{k}\},\quad \text{if $n$ is  odd}.\end{matrix}\right.$$
\end{theorem}

{\bf Proof.} Let $\g\subset\co(1,n+1)$ be the holonomy algebra of $(M,c,\nabla)$. It is clear that the Lie algebra $\co(1,n+1)$ does not annihilate any non-zero element of the module $\Delta_{1,n+1}$ for arbitrary~$w$. 

Let $\g$ be the first algebra from Theorem \ref{class2Th}.
Suppose that it annihilates a non-zero element in $\Delta_{1,n+1}$. Then $\so(1,k+1)\oplus\so(n-k)$ annihilates this element as well. 
On the other hand, depending on $n$ and $k$,  $\Delta_{1,n+1}$ considered as the  
$\so(1,k+1)\oplus\so(n-k)$-module is either isomorphic to $\Delta_{1,k+1}\otimes\Delta_{n-k}$ or to the direct sum of two copies of this module. This gives a contradiction since due to the restriction on $n$ and $k$ it is not possible that both $\so(1,k+1)$ annihilates a non-zero element of $\Delta_{1,k+1}$ and $\so(n-k)$ annihilates a non-zero element of $\Delta_{n-k}$.

Now we consider holonomy algebras from Theorem \ref{class3Th}. Let $\g$ be one of these algebras.

Under the identification $\mathfrak{so}(1,n+1)\cong
\Lambda^2\mathbb{R}^{1,n+1}$ given by \eqref{isom_so}, an element
$(0,a,A,X)\in\mathfrak{so}(1,n+1)_{\Real p}$ corresponds to the bivector
$$-ap\wedge q+A-p\wedge X.$$ We follow the computations from \cite{L01}. It holds $$\Delta_{1,n+1}\cong
\Delta_{n} \otimes \Delta_{1,1},\quad \Delta_{1,1}\cong
\mathbb{C}^2.$$ Consider the vectors $e_-=\frac{\sqrt{2}}{2}(p-q)$
and $e_+=\frac{\sqrt{2}}{2}(p+q)$ and the orthonormal basis $e_-
,e_+,e_1,...,e_n$ of $\mathbb{R}^{1,n+1}$. Note that
$p=\frac{\sqrt{2}}{2}(e_-+e_+)$ and
$q=\frac{\sqrt{2}}{2}(e_--e_+)$.

Let $\psi\in \Delta_{1,n+1}$ be a spinor annihilated by $\g$. We
may write
$$\psi=\psi_+\otimes u(1)+\psi_-\otimes u(-1),$$ where
$\psi_\pm\in\Delta_{n}$. Using \eqref{Phi1}, \eqref{Phi2} and the
computations from \cite{L01} it is easy to get that $$
(e_1\wedge p)\cdot \psi=\frac{\sqrt{2}}{2}e_1\cdot(e_-+e_+) \cdot
\psi=\sqrt{2}(e_1\cdot \psi_-)\otimes u(1).$$ Hence the equality
$(e_1\wedge p)\cdot \psi=0$ imply $e_1\cdot \psi_-=0$.
Since non-zero vectors from $\Real^n$ act in $\Delta_n$ as  isomorphisms, we get
$\psi_-=0$. Thus, $\psi=\psi_+\otimes u(1)$, i.e., $\psi\in\Delta_{n}\otimes u(1)$. Moreover, it is clear that   $\Delta_{n}\otimes u(1)$ is a trivial
$p\wedge\mathbb{R}^n$-module.  Next, $$(p\wedge
q)\cdot(\psi_+\otimes u(1))=2\psi_+\otimes u(1),\quad A(\psi_+\otimes
u(1))=A(\psi_+)\otimes u(1)$$ for all $\psi_+\in\Delta_n$ and
$A\in\mathfrak{so}(n)$.

Suppose that $\Real\id\subset\g$, i.e., $\g=\Real\id\oplus\tilde\g$, where $\tilde \g\subset\so(1,n+1)_{\Real p}$ is the holonomy algebra of a Lorentzian manifold. If $\g$ annihilates a spinor $\psi$, then both $\id$ and $\tilde\g$ annihilate $\psi$. The first condition implies that $w=0$.  According to \cite{L01,L07}, the second condition implies that $\tilde\g=\h\zr\Real^n$, where $\h\subset\so(n)$ 
is the holonomy algebra of a Riemannian manifold carrying a parallel spinor. We see that $\g$ is the first algebra from the statement of the theorem with $w=0$.

Consider a holonomy algebra from Theorem \ref{class3Th} defined by a map $\theta$.  
Suppose that it annihilates a non-zero spinor $\psi$ of weight $w$. Let $(\theta(A),a,A,0)\in\g$. Then 
$\psi=\psi_+\otimes u(1)$ and
$$A\psi_+=(2a-w\theta(A))\psi_+.$$
The element $A\in\so(n)$ may be written in the canonical from
$$A=\sum_{k=1}^mc_ke_{2k-1}\wedge e_{2k},$$ where $2m\leq n$, $c_k\in\Real$.
Suppose that $n$ is even.
According to Section \ref{sec1.1}, the element $e_{2k-1}\wedge e_{2k}$ acts on $\Delta_n$ as the endomorphism
$$\Phi_{n}(e_{2k-1})\Phi_{n}(e_{2k})= E\otimes\cdots\otimes E\otimes
UV\otimes \underbrace{E\otimes \cdots\otimes
	E}_{(k-1)-times},\quad 
UV=\left(\begin{array}{cc}0&-1\\1&0\end{array}\right).$$
This implies that $A$ acting on $\Delta_n$ has no non-zero real eigenvalues. The same can be shown if  $n$ is odd. We get  $$w\theta(A)=2a,\quad A\psi_+=0.$$ for all $A\in\h$.  This implies that $\h$ annihilate $\psi_+$.
As we have seen in Section \ref{sec2}, the center of $\h$ is trivial.
This implies that $\theta=0$. Thus, $\g$ is $\g^{\beta,\theta,1,\h}$ with $\theta=0$, and $\beta=\frac{2}{w}$. 

Thus  the algebras of Theorem \ref{class3Th} give us the first algebra from the statement of the theorem.  This algebra annihilates a spinor $\psi$ as above if and only if the corresponding subalgebra $\h$ annihilates the spinor $\psi_+$. 

Let $\g$ be the second algebra from Theorem \ref{class2Th}.
Now it is clear that if this algebra annihilates a non-trivial element in  $\Delta_{1,n+1}$ then $\so(n-k)$ must annihilate a spinor in $\Delta_{n-k}$. This is possible only if $k=n-1$. We obtain the algebra $\g^\mathfrak{k}$. 
The third algebra from Theorem \ref{class2Th} satisfies $\g=\Real\id\oplus\tilde\g$, where $\tilde \g\subset\so(1,n+1)_{\Real p}$. The above arguments show that $\g$ does not annihilate any non-zero spinor.

Finally, the 
statement about the dimension of parallel spinors for the algebra $\g^\mathfrak{k}$  follows from the obvious fact that $\Delta_n$ considered as the $\so(n-1)$-module coincides with $\Delta_{n-1}$ if $n$ is even, and   $\Delta_n$ is the direct sum of two copies of $\Delta_{n-1}$ if $n$ is odd. The theorem is proved. \qed

Note that if the holonomy algebra of a simply connected Weyl manifold $(M,c,\nabla)$ is $\g^\mathfrak{k}$, then $(M,c,\nabla)$ admits a parallel null vector field; $(M,c,\nabla)$ is a conformal product in the sense of \cite{Moroianu1} with a 1-dimensional factor, and $(M,c,\nabla)$ admits also a recurrent space-like  vector field.

Theorem \ref{thHolSpinor} implies that if a simply connected  manifold $(M,c,\nabla)$ with a non-closed Weyl structure admits a non-zero weighted parallel spinor, then it admits a non-zero  recurrent null vector field. This may be seen in the following way. The spinor module $\Delta_{1,n+1}$ admits an $\so(1,n+1)$-invariant Hermitian product $b$. A spinor $\psi$ defines the Dirac current $V_\psi$ by the equality $$g(V_\psi,X)=-b(X\cdot\psi,\psi).$$ 
In \cite{MOP}, it is shown that if $\psi$ is a weighted parallel spinor, then $V_\psi$ is a recurrent vector field. Next, it is well-known (see, e.g., \cite{MOP} or \cite{L01}) that in the Lorentzian signature it holds $g(V_\psi,V_\psi)\leq 0$, and $V_\psi(x)=0$ ($x\in M$) if and only if $\psi(x)=0$. Since $\psi$ is parallel, it is nowhere vanishing,  and hence the vector field $V_\psi$ is  either   null, or it is time-like. The second case would imply that the holonomy algebra $\g$ of $(M,c,\nabla)$ preserves a time-like line in $\Real^{1,n+1}$, and by Theorem \ref{class2Th}, $\g$ would coincide with $\Real\id_{\Real^{1,n+1}}\oplus\so(n+1)$, but since $n\geq 1$, in that case $(M,c,\nabla)$ does not admit any non-zero   weighted parallel spinor. Thus, the Dirac current is null, and consequently it is proportional to the vector field  $\partial_v$.

\section{Kundt and Walker structures for Weyl manifolds}\label{sec4}

Our next task is to describe the structures of the spaces with the holonomy algebras from Theorem \ref{thHolSpinor}. For that we introduce in this section  special local coordinates for Weyl manifolds admitting parallel distributions of null lines.

On the space $\Real^{n+2}$ consider the coordinates $v,x^1,\dots,x^n,u$ and the metric
\begin{equation}\label{metric} g=2dvdu+h+2 Adu +H (du)^2,\end{equation}
where $$h=
h_{ij}dx^idx^j,\quad\partial_vh_{ij}=0,$$ is a $u$-family of Riemannian metrics, and $$A=A_idx^i$$ is a 1-form. Such metrics $g$ were introduced in \cite{Kundt} and are called Kundt metrics. If the functions $A_i$ are independent of the coordinate $v$, then $g$ is a Walker metric \cite{Walkerbook}.
In that case, the null vector field $\partial_v$ is recurrent with respect to the Levi-Civita connection, and it generates a parallel distribution of null lines. Each Lorentzian manifold with a parallel distribution of null lines is locally given  by such  metric.

Suppose now that $(M,c,\nabla)$ is a Lorentzian Weyl manifold that admits a parallel distribution of null lines, i.e., the holonomy algebra of $(M,c,\nabla)$ is contained in $\co(1,n+1)_{\Real p}$. The local form of 
general $(M,c,\nabla)$ is given in the following theorem.

\begin{theorem}\label{localform1}
Let $(M,c,\nabla)$ be a Lorentzian Weyl manifold of dimension $n+2$  admitting a parallel distribution of null lines. Then  around each point of $M$ there exist local coordinates $v,x^1,\dots,x^n,u$ and a Kundt metric $g\in c$  
 given by  \eqref{metric} such that the corresponding 1-form $\omega$ satisfies $\omega(\partial_v)=0.$
\end{theorem}

If $(M,c,\nabla)$ satisfies an additional condition on the curvature, then we obtain a finer description of the local structure of $(M,c,\nabla)$.

\begin{theorem}\label{localform2}
Let $(M,c,\nabla)$ be a Lorentzian Weyl manifold of dimension $n+2$  admitting a parallel distribution of null lines $\ell$. Suppose that the curvature tensor $R$ of the connection $\nabla$ satisfies the condition
\begin{equation}\label{conditionR} R(X,Y)p=0, \quad \forall\, X,Y\in\ell^\bot,\,p\in\ell.\end{equation}

Then  around each point of $M$ there exist local coordinates $v,x^1,\dots,x^n,u$ and a Walker metric $g\in c$   (given by \eqref{metric} with $\partial_v A_i=0$) such that $\omega=fdu$ for a function~$f$. Moreover, the coordinates may be chosen in such a way that~$A=0$.
	\end{theorem}

The condition \eqref{conditionR} is satisfied for Weyl structures with a broad class of holonomy algebras:    

\begin{theorem}\label{th_usl_dlja} 
	Let $(M,c,\nabla)$ be a Lorentzian Weyl manifold of dimension $n+2$  admitting a parallel distribution of null lines $\ell$. Let $\g\subset \co(1,n+1)_{\Real p}$ be its holonomy algebra. If the projection of  $\g$ to the subalgebra 
	$\so(n)\subset\co(1,n+1)_{\Real p}$ is different from $\so(n)$, then the curvature tensor $R$ of $\nabla$ satisfies the condition \eqref{conditionR}, and, in particular, by Theorem \ref{localform2}, there is a Walker metric in the conformal class $c$.  
\end{theorem}

{\bf Proof of Theorems \ref{localform1} and \ref{localform2}.} Let $\ell$ be a parallel distribution of null lines  on $(M,c,\nabla)$.
It is clear that the orthogonal complement $\ell^\bot$ (taking with respect an arbitrary $g\in c$) is also parallel. Since the connection $\nabla$ is tosion-free, the distribution $\ell^\bot$ is involutive. The flag of involutive distributions $\ell\subset\ell^\bot$ defines a flag of foliations. Let $v,x^1,\dots,x^n,u$ be local coordinates corresponding to that flag of foliations. This means that the vector field $\partial_v$ spans the distribution $\ell$, while the vector fields 
$\partial_v,\partial_{x^1},\dots,\partial_{x^n}$ span the distribution $\ell^\bot$. Let $g\in c$ be a metric such that $g(\partial_v,\partial_u)=1$. Then it is given by \eqref{metric} with $h$ possibly depending on $v$.

Applying the vector field $\partial_v$ to the equality $1=g(\partial_v,\partial_u)$ and using \eqref{nablag},
we get $$0=2\omega(\partial_v)g(\partial_v,\partial_u)+g(\nabla_{\partial_v}\partial_v,\partial_u)+g(\partial_v,\nabla_{\partial_v}\partial_u).$$
This implies
\begin{equation}\label{eq1}
2\omega(\partial_v)+\Gamma^v_{vv}=0.\end{equation}
In the same way we obtain
\begin{equation}\label{eq2}\partial_vh_{ij}=2\omega(\partial_v)h_{ij},\quad  \partial_vA_i=2\omega(\partial_v)A_i+\Gamma^v_{vi},\quad 2\omega(\partial_i)+\Gamma^v_{vi}=0.\end{equation}

 The arbitrary coordinate transformation that preserves the above flag of foliations is of the form 
\begin{equation}\label{coordtransf}\tilde v=\tilde v(v,x^1,\dots,x^n,u),\quad 
 x^{\tilde i}=\tilde x^{\tilde i}(x^1,\dots,x^n,u),\quad \tilde u=\tilde u(u).\end{equation}

The following relations for the Christoffel symbols hold:
$$\Gamma^v_{vv}=\Gamma^{\tilde v}_{{\tilde v}{\tilde v}}\partial_v{\tilde v}+\partial_v\ln|\partial_v\tilde v|,$$
$$\Gamma^v_{vi}=\Gamma^{\tilde v}_{{\tilde v}{\tilde v}}\partial_{i}{\tilde v}+\Gamma^{\tilde v}_{{\tilde v}{\tilde k}}\partial_{i}{ x^{\tilde k}}+\partial_{x^i}\ln|\partial_v\tilde v|.$$

The condition $$\Gamma^{\tilde v}_{{\tilde v}{\tilde v}}=0$$ is equivalent to the differential equation $$\partial_v \ln|\partial_v{\tilde v}|=\Gamma^v_{vv},$$
which is integrable. Hence there exist local coordinates $v,x^1,\dots,x^n,u$ such that $\Gamma^v_{vv}=0$. Choosing the new metric $g$ in such a way that again $g(\partial_v,\partial_u)=1$ we get the proof of  Theorem \ref{localform1}.

The conditions $$\Gamma^{\tilde v}_{{\tilde v}{\tilde v}}=\Gamma^{\tilde v}_{{\tilde v}{\tilde i}}=0$$ are equivalent to the system of partial differential equations
$$\partial_v \ln|\partial_v{\tilde v}|=\Gamma^v_{vv},\quad \partial_{i} \ln|\partial_v{\tilde v}|=\Gamma^v_{vi}.$$
The integrability conditions for this system are 
$$\partial_v\Gamma^v_{vi}=\partial_i\Gamma^v_{vv}, \quad \partial_j\Gamma^v_{vi}=\partial_i\Gamma^v_{vj}.$$
We claim that these conditions are exactly the condition \eqref{conditionR}.
We use the convention $$R(\partial_{a},\partial_{b})\partial_{c}=R^d_{cab}\partial_{d},$$ where 
the indices take values $v,1,\dots,n,u$. Since the vector field $\partial_v$ is recurrent, the condition~\eqref{conditionR} takes the form
$$R^v_{vvi}=R^v_{vij}=0.$$
It is easy to check that 
$$R^v_{vvi}= \partial_v\Gamma^v_{vi}-\partial_i\Gamma^v_{vv},\quad
R^v_{vij}= \partial_i\Gamma^v_{vj}-\partial_j\Gamma^v_{vi}.$$
This proves the claim.  
Thus if the condition \eqref{conditionR} is satisfied, then there exist coordinates $v,x^1,\dots,x^n,u$ such that 
$\Gamma^v_{vv}=\Gamma^v_{vi}=0$. Choosing again the new metric $g$ in such a way that  $g(\partial_v,\partial_u)=1$, and using \eqref{eq1}, \eqref{eq2},  we obtain the proof of Theorem \ref{localform2} (the last statement of Theorem~\ref{localform2} follows from the results of~\cite{GL10}). \qed

{\bf Proof of Theorem \ref{th_usl_dlja}} immediately follows from \cite[Th.~8]{Andr}, where it is trivial to see that the condition on the holonomy algebra from the statement of Theorem \ref{th_usl_dlja} implies \eqref{conditionR}. \qed

\section{Local structure of Lorentzian Weyl manifolds with parallel spinors}\label{sec5}

In this section we describe the local form of the Weyl structures with the holonomy algebras from Theorem \ref{thHolSpinor}, i.e., we describe the Weyl structures admitting weighted parallel spinors.

For a metric \begin{equation}\label{WmetricA=0}
g=2dvdu+h+H(du)^2,\quad h=h_{ij}dx^idx^j,\quad \partial_v h_{ij}=0,\end{equation}
and a form $\omega=fdu$ consider the corresponding Weyl connection $\nabla$ and the flag of parallel distributions $\ell\subset\ell^\bot$. Consider the quotient ${\mathcal E}=\ell^\bot/\ell$.
The connection $\nabla$ induces the connection $\nabla^{\mathcal E}$ on the bundle~${\mathcal E}$.  Likewise, the connection $\nabla^g$ induces the connection $\nabla^{g,\mathcal E}$ on ${\mathcal E}$. Consider the vector fields $\partial_i$ as the sections of ${\mathcal E}$. Let $\Gamma^i_{aj}$ and $\bar\Gamma^i_{aj}$ be the Christoffel symbols of the connections $\nabla^{\mathcal E}$ and $\nabla^{g,\mathcal E}$, respectively. It holds 
\begin{equation}
\label{GammaE}
\Gamma^i_{vj}=\bar\Gamma^i_{vj}=0,\quad \Gamma^i_{kj}=\bar\Gamma^i_{kj},\quad \Gamma^i_{uj}=\bar\Gamma^i_{uj}+f\delta^i_j.
\end{equation}
This shows that the $\so(n)$-projections of the holonomy algebras of the connections $\nabla$, $\nabla^g$, $\nabla^{\mathcal E}$  coincide with the holonomy algebra of the connection $\nabla^{g,\mathcal E}$. Now, the $\so(n)$-projection of the holonomy algebra of the metric $g$ may be found using the algorithms from~\cite{GalLMLP}.

\begin{theorem}\label{thconnect1}
Let $(M,c,\nabla)$ be a simply connected manifold with a non-closed Weyl spin structure of Lorentzian signature. Then $(M,c,\nabla)$ admits a non-zero weighted parallel spinor of weight~$w$ if and only if the holonomy algebra $\h\subset\so(n)$ of the connection $\nabla^{g,\mathcal E}$ is the holonomy algebra of a Riemannian spin manifold admitting a non-zero parallel spinor,  around each point of $M$ there exist local coordinates $v,x^1,\dots,x^n,u$ and  a metric $g\in c$  such that 
$$g=2dvdu+h+H(du)^2,$$
where $H$ is a function, $h=h_{ij}dx^idx^j$, 
 $\partial_vh_{ij}=0$, and the corresponding 1-form $\omega$ satisfies \begin{equation}\label{condtheta1}\omega=fdu,\quad (2+w)f=\partial_vH.\end{equation} 
\end{theorem}

{\bf Proof.}  
Suppose that $(M,c,\nabla)$ admits a non-zero weighted parallel spinor.  
 By Theorems \ref{thHolSpinor}, \ref{th_usl_dlja} and \ref{localform2}, around each point of $M$ there exist coordinates $v,x^1,\dots,x^n,u$ and a metric $g\in c$ such that $g$ is a Walker metric with $A=0$, and $\omega=fdu$.  Consider the frame
$$p=\partial_v,\quad X_i=\partial_i, \quad q=\partial_u-\frac{1}{2}H\partial_v.$$
Let $E$ be the distribution spanned by the vector fields $X_1,\dots,X_n$.
By Theorem \ref{thHolSpinor}, the holonomy algebra $\g$ of $\nabla$ is either $\g^{w,\h}$ or $\g^{\mathfrak{k}}$. In both cases the projection of $\g$ to $\Real\id\oplus\Real(0,1,0,0)\subset\co(1,n+1)_{\Real p}$    is equal to $\Real(2,w,0,0)$.
The form of the holonomy algebra remains unchanged under any transformation of the Witt basis from Section \ref{sec2}. This and the Ambrose-Singer theorem imply 
that the condition on the projection of the holonomy algebra to
$\Real\id\oplus\Real(0,1,0,0)\subset\co(1,n+1)_{\Real p}$  is equivalent to
 the equalities
\begin{equation}\label{uslovijaR}g(R(X,Y)p,q)g(V,V)=(2+w)g(R(X,Y)V,V),\,\,
\forall X,Y\in \Gamma(TM),\,\, V\in\Gamma(E) . \end{equation}
These equalities are equivalent to
$$R^v_{vab}=(2+w)R^i_{iab},\quad a,b=v,1,\dots,n,u,\quad\text{(no summation over $i$)}.$$
The simple direct computations give the expressions
$$R^v_{vvj}=R^i_{ivj}=R^v_{vjk}=R^i_{ijk}=0,$$
$$R^v_{vvu}=\frac{1}{2}\partial^2_vH,\quad R^i_{ivu}=\frac{1}{2}\partial_vf,\quad 
R^v_{vju}=\frac{1}{2}\partial_v\partial_jH,\quad R^i_{iju}=\frac{1}{2}\partial_jf.$$
This implies the equality 
$(2+w)f=\partial_vH+F(u)$, where $F(u)$ is a function. 
If $w\neq -2$, then changing $\omega$ to $\omega -d(G(u))$ for a proper $G(u)$, we may assume that $F(u)=0$. Similarly, if $w=-2$, then $H=-vF(u)+H_0$, $\partial_vH_0=0$. The function $F(u)$ does not influence the curvature of the metric $g$ \cite{UMN}, consequently, the coordinates may be choosen in such a way that  $F(u)=0$.

Let us prove the inverse implication. It is clear that the vector field $\partial_v$ is recurrent. Hence the holonomy algebra $\g$ is either $\g^{\mathfrak{k}}$ or the third algebra from Theorem \ref{class2Th}, or it is one of the algebras from Theorem~\ref{class3Th}. In the first case there is nothing to prove. Consider the other cases. Since $\h\subset\so(n)$ is the holonomy algebra of a Riemannian manifold carrying a non-zero parallel spinor, the center of $\h$ is trivial. This implies that $\g$ is one of the algebras $\g^{\Real,1,\h}$, $\g^{\Real,2,\h}$, $\g^{\theta,\beta,1,\h}$ with $\theta=0$, or it is the third algebra from Theorem \ref{class2Th}. The conditions on $f$ and $H$, and the proof of the first implication imply that $\g$ is $\g^{w,\h}$.
 This proves the theorem.  \qed

In the case of the second holonomy algebra from Theorem~\ref{thHolSpinor} we may say more.

\begin{theorem}\label{thconnect2}
	Let $(M,c,\nabla)$ be a simply connected  manifold with a non-closed Weyl structure.
		 Then its holonomy algebra equals $\g^{\mathfrak{k}}$ with $\mathfrak{k}\subset\so(n-1)$ being  a Riemannian holonomy algebra if and only if around each point of $M$  there exist local coordinates $v,x^1,\dots,x^n,u$ and a metric $g\in c$   such that 
	$$g=2dvdu+h+H(du)^2,$$
	where $H=H(x^1,\dots,x^{n-1},u)$ is a function, $$h=\sum_{i,j=1}^{n-1}h_{ij}(x^1,\dots,x^{n-1},u)dx^idx^j+e^{-2F}(dx^n)^2,$$
	$$\omega =fdu,\quad f=\partial_uF,$$
	and $F=F(x^n,u)$ is a function. 
	 \end{theorem}

{\bf Proof.} Suppose that the holonomy algebra of a Weyl structure is  $\g^\mathfrak{k}$. It is clear that there exist a metric $g$ and a form $\omega$ as in Theorem~\ref{thconnect1}.  Since $\mathfrak{k}\subset\so(n-1)\subset  
\so(n)$ preserves the vector subspace $\Real^{n-1}\subset\Real^n$, according to the discussion before Theorem \ref{thconnect1} and to  results from~\cite{Boubel}, the cooridnates may be chosen in such a way that
$$h=\sum_{i,j=1}^{n-1}h_{ij}(x^1,\dots,x^{n-1},u)dx^idx^j+h_{nn}(x^n,u)(dx^n)^2.$$ 
It holds 
$$\Gamma^v_{vn}=0,\quad \Gamma^v_{nn}=-\frac{1}{2}\partial_nH,\quad \Gamma^v_{in}=0,\,\ 1\leq i\leq n-1,\quad 
\Gamma^v_{un}=- \frac{1}{2}\partial_u h_{nn}-fh_{nn}.$$
These equalities and condition that the vector field $\partial_n$ is recurrent imply the proof of the first implication. The proof of the inverse implication is obvious. \qed 

Let us consider examples of Weyl manifolds with weighted parallel spinors. 

\begin{ex} Let $(N,h)$ be a simply connected Riemannian spin  manifold of dimension $n-1$ carrying a non-zero parallel spinor. Let $F(x^n,u)$ be a function such that $\partial_{x^n}\partial_vF\neq 0$. 	
	Let $$M=\Real\times N\times\Real\times \Real,$$
	$$g=2dvdu+h+e^{-2F}(dx^n)^2,\quad c=[g],$$
	$$\omega=\partial_uFdu,$$
	and $\nabla$ be the Weyl connection defined by $g$ and $\omega$ as in \eqref{formulaK}. Then $(M,c,\nabla)$ is a non-closed Weyl structure with the holonomy algebra $\g^{\mathfrak{k}}$, where $\mathfrak{k}\subset\so(n-1)$ is the holonomy algebra of $(N,h)$, and consequently  $(M,c,\nabla)$ carries a non-zero parallel spinor of weight $-2$. The dimension of the space of weighted parallel spinors of weight $-2$ is the dimension of the space of parallel spinors on $(N,h)$ if $n$ is even, and it is two times so big if $n$ is odd.
	\end{ex}

\begin{ex} Let $(N,h)$ be a simply connected Riemannian spin manifold of dimension $n$ carrying a non-zero parallel spinor. Suppose that $w\neq -2$. Let $H$ be a function with $\partial^2_v H\neq 0$, and 
	$\omega =fdu$, where $f=\frac{1}{2+w}\partial_vH$.
	Let $$M=\Real\times N\times \Real,$$
	$$g=2dvdu+h+2H(du)^2,\quad c=[g],$$
and $\nabla$ be the Weyl connection defined by $g$ and  $\omega$ as in \eqref{formulaK}. Then $(M,c,\nabla)$ is a non-closed Weyl structure with the holonomy algebra $\g^{w,\mathfrak{h}}$, where $\mathfrak{h}\subset\so(n)$ is the holonomy algebra of $(N,h)$, and consequently  $(M,c,\nabla)$ carries a non-zero parallel spinor of weight $w$. The dimension of the space of weighted parallel spinors of weight $w$ equals the dimension of the space of parallel spinors on~$(N,h)$.
 \end{ex}

\section{Einstein-Weyl equation}\label{secEW}

Let $(M,c,\nabla)$ be a Weyl manifold. We will say that it is Einstein-Weyl if it holds
$$\Ric^s=\Lambda g,$$
where $\Ric^s$ is  the symmetric part of the Ricci tensor of the connection $\nabla$, $g\in c$, and $\Lambda$ is a function on $M$.

For a metric \eqref{WmetricA=0} and a local function $H$ on $M$ let
$\Delta H=h^{ij}\nabla_i\nabla_jH$ denote the Laplace-Beltrami operator with respect to the $u$-family of the Riemannian metrics $h$ applied to $H$. 
A dot over a function will denote the partial derivative with respect to the coordinate $u$.

\begin{theorem}\label{thEWLam}
	Let $(M,c,\nabla)$ be a simply connected manifold with a non-closed  Weyl spin structure of Lorentzian signature $(1,n+1)$, $n\ge1 $. Then $(M,c,\nabla)$ is Einstein-Weyl and it admits a non-zero weighted parallel spinor of weight $w$ if and only if $$w=n-2,$$
		 the holonomy algebra $\h\subset\so(n)$ of the connection $\nabla^{g,\mathcal E}$ is the holonomy algebra of a  Riemannian spin manifold admitting a non-zero parallel spinor,  around each point of $M$ there exist local coordinates $v,x^1,\dots,x^n,u$ and   a metric $g\in c$   such that 
	$$g=2dvdu+h+H(du)^2,$$
	where   $h=h_{ij}dx^idx^j$, 
	$\partial_vh_{ij}=0$,  the corresponding 1-form $\omega$ satisfies \begin{equation*}\omega=\frac{1}{n}\partial_v Hdu,
		\end{equation*} 
		and the following equation holds:
		\begin{equation}\label{Ricuu=0}
		\frac{2(n-2)}{n}(\partial^2_v H)-2H\partial^2_vH+4\partial_u\partial_v H+2\Delta H-h^{ij}h^{kl}\dot h_{ik}\dot h_{jl}+2h^{ij} \ddot h_{ij}+h^{ij}\dot h_{ij}\partial_vH=0.
		\end{equation}
						\end{theorem}

{\bf Proof.} Let us consider a Weyl structure as in Theorem \ref{thconnect1} and impose the Einstein-Weyl equation.
The symmetric part $\Ric^s$ of the connection $\nabla$ may be written in the form
$$\Ric^s=\overline \Ric+\widehat\Ric,$$
where $\overline\Ric$ is the Ricci tensor of the metric $g$. The assumption on the subalgebra $\h\subset\so(n)$ and the  computations from \cite{UMN} imply that
\begin{align*} \overline\Ric_{vv}&=\overline\Ric_{vi}=\overline\Ric_{ij}=0,\\
\overline\Ric_{vu}&=-\frac{1}{2}\partial^2_vH,\quad \overline\Ric_{iu}=-\frac{1}{2}\partial_v\partial_iH,\\
\overline\Ric_{uu}&=-\frac{1}{2}H\partial^2_vH+\frac{1}{2}\Delta H-\frac{1}{4}h^{ij}h^{kl}\dot h_{ik}\dot h_{jl}+\frac{1}{2}h^{ij} \ddot h_{ij}+\frac{1}{4}h^{ij}\dot h_{ij}\partial_vH.
\end{align*}
It is easy to check directly that 
\begin{align*}
\widehat\Ric_{vv}&=\widehat\Ric_{vi}=0,\quad \widehat\Ric_{ij}=h_{ij}\partial_v f,\\
\widehat \Ric_{vu}&=\frac{n+2}{2}\partial_v f,\quad \widehat \Ric_{iu}=\frac{n}{2}\partial_i f,\\
 \widehat \Ric_{uu}&=n\dot f-nf^2+H\partial_vf+\frac{n}{2}f\partial_v H.
\end{align*}
Now, the equation $\Ric^s_{ij}=\Lambda h_{ij}$ is equivalent to $\Lambda=\partial_v f$. Using this, we see that
the equations $\Ric^s_{vu}=\Lambda$ and  $\Ric^s_{iu}=0$ are equivalent to
$$n\partial_vf-\partial^2_vH=0,\quad n\partial_i f-\partial_i\partial_vH=0.$$
Recall that it holds $(w+2)f=\partial_v H$. We obtain the equations
$$(w+2-n)\partial_vf=(w+2-n)\partial_i f=0.$$
Since $d\omega\neq 0$, we get $w=n-2$.
The rest of the proof of the theorem follows directly from the equality $\Ric^s_{uu}=0$. In fact, in \cite{MOP} it is shown that the existence of a non-zero weighted parallel spinor of weight $w=n-2$ implies that the tensor $\Ric-\Lambda g$ is totally isotropic ($\Lambda$ is the scalar curvature of $\nabla$ devided by $n+2$), this condition is equivalent to the equality  $\Ric^s_{uu}=0$.   \qed  

\begin{cor}\label{corEW1}
	Let $(M,c,\nabla)$ be a simply connected manifold with a non-closed  Weyl spin structure of Lorentzian signature $(1,n+1)$, $n\ge1 $. If $(M,c,\nabla)$ is Einstein-Weyl and  admits a non-zero weighted parallel spinor of weight $w$, then $w=n-2$, and the holonomy algebra of $(M,c,\nabla)$ is $\g^{w,\h}$, where $\h\subset\so(n)$ is the holonomy algebra of a Riemannian spin manifold admitting a non-zero parallel spinor.
	\end{cor}

\begin{cor}\label{corEW0}
	Let $(M,c,\nabla)$ be a simply connected manifold with a non-closed  Weyl spin structure of Lorentzian signature $(1,n+1)$, $n\ge1 $. Then $(M,c,\nabla)$ is Einstein-Weyl with zero scalar curvature and  admits a non-zero weighted parallel spinor of weight $w$ if and only if $w=n-2$,
	the holonomy algebra $\h\subset\so(n)$ of the connection $\nabla^{g,\mathcal E}$ is the holonomy algebra of a  Riemannian spin manifold admitting a non-zero parallel spinor,  around each point of $M$ there exist local coordinates $v,x^1,\dots,x^n,u$,   a metric $g\in c$ and a function $f$ with $\partial_vf=0$   such that 
	$$g=2dvdu+h+H(du)^2,$$
	where $$H=nfv+H_0,\quad \partial_vH_0=0,$$  $h=h_{ij}dx^idx^j$, 
	$\partial_vh_{ij}=0$,  the corresponding 1-form $\omega$ satisfies \begin{equation*}\omega=fdu,
	\end{equation*} 
	and the following equations hold:
	\begin{align*}\Delta f&=0,\\ 
	2n(n-2)f^2-4n\dot f+2\Delta H_0-h^{ij}h^{kl}\dot h_{ik}\dot h_{jl}+2h^{ij} \ddot h_{ij}+nfh^{ij}\dot h_{ij}&=0.
	\end{align*}
\end{cor}

\begin{ex} Let $(N,h)$ be a simply connected Riemannian spin manifold of dimension $n$ carrying a non-zero parallel spinor and with a non-constant harmonic function $F$. Let $w=n-2$.  If $n=2$, then we set $f=F$. If $n\neq 2$, then we set $f=-\frac{2}{(n-2)u}F.$
		Let  
	$\omega =fdu$.
	Let $$M=\Real\times N\times \Real_+,$$
	$$g=2dvdu+h+2nf(du)^2,\quad c=[g],$$
	and $\nabla$ be the Weyl connection defined by $g$ and  $\omega$ as in \eqref{formulaK}. Then $(M,c,\nabla)$ is an Einstein-Weyl non-closed Weyl structure with zero scalar curvature and the holonomy algebra $\g^{w,\mathfrak{h}}$, where $\mathfrak{h}\subset\so(n)$ is the holonomy algebra of $(N,h)$, and   $(M,c,\nabla)$ curries a non-zero parallel spinor of weight $w=n-2$. The dimension of the space of weighted parallel spinors of weight $w=n-2$ is equal to the dimension of the space of parallel spinors on $(N,h)$.
\end{ex} 

\begin{ex} 
Let $(M,c,\nabla)$ be a simply connected manifold with a non-closed  Weyl spin structure of Lorentzian signature $(1,n+1)$, $n\ge1 $.
Suppose that $(M,c,\nabla)$  admits a non-zero weighted parallel spinor of weight $w$. Suppose that $1\leq n\leq 3$.  The above condition on $\h\subset\so(n)$ implies that $\h=0$, and $h_{ij}=\delta_{ij}$. The Einstein-Weyl equation implies that
$$g=2dvdu+\delta_{ij}dx^idx^j+H(du)^2,$$
$$\omega=\frac{1}{n}\partial_v Hdu,$$
and $H$ satisfies  
$$\frac{2(n-2)}{n}(\partial^2_v H)-2H\partial^2_vH+4\partial_u\partial_v H+2\delta^{ij}\partial_i\partial_j H=0.$$
\end{ex}
The  
Einstein-Weyl equation in dimension 3 ($n=1$) and its relation to various geometric structures have been carefully studied in many works, for example see \cite{DGS,DFK,DMT}. In that case, the equation on $H$ is equivalent to the dispersionless Kadomtsev-Petviashvili equation. 
The case $n=2$ is considered in \cite{MOP}.
Note that if $n=4$, then either $\h=0$, or $\h=\mathfrak{su}(2)$, this case was studied in \cite{MOP}.

\medskip

In conclusion let us compare the results of the current papers with the results from \cite{MOP}. In \cite{MOP}, were considered Einstein-Weyl manifolds $M$ with parallel spinors of weight $\dim M-4$. Using the  the parallel spinor equation, it was shown that the conformal class contains a Walker metric. The Einstein-Weyl equation was carefully studied in dimensions 4 and 6. In the present paper we considered parallel spinors of arbitrary weight showing that the Einstein-Weyl case corresponds to spinors of weight $\dim M-4$. We found the holonomy algebras and the local form for all Weyl structures carrying parallel spinors of arbitrary weight. We introduced the Kundt and Walker coordinates for all Weyl connections admitting parallel distributions of isotropic lines (not necessary carrying a parallel spinor).

\vskip0.5cm

{\bf Data Availability Statement.} 
Data sharing not applicable to this article as no datasets were generated or analysed during the current study.

\vskip0.5cm

{\bf Acknowledgements.} The authors are thankful to Maciej Dunajski for useful email communications and to Ioannis Chrysikos for useful suggestions. The authors are grateful to the anonymous referee for careful reading of the paper and valuable comments that greatly affected the final appearance of the paper. A.D. was supported by  grant MUNI/A/1160/2020 of Masaryk University.   A.G. was supported by grant no. 18-00496S of the Czech Science Foundation.

\end{document}